\documentclass{amsart}
\usepackage{amssymb}
\usepackage{epsfig,pb-diagram}

\title{Growth Diagrams for the Schubert  Multiplication}

\author{Cristian Lenart}

\address{Department of Mathematics and Statistics, State University of New York at Albany, Albany, NY 12222}
\email{lenart@albany.edu}

\thanks{ C. L. was partially supported by National Science Foundation 
grant DMS-0701044}

\newlength{\cellsize}
\newcommand\tableau[1]{
\vcenter{
\let\\=\cr
\baselineskip=-16000pt
\lineskiplimit=16000pt
\lineskip=0pt
\halign{&\tableaucell{##}\cr#1\crcr}}}

\newcommand{\tableaucell}[1]{{%
\def \arg{#1}\def \void{}%
\ifx \void \arg
\vbox to \cellsize{\vfil \hrule width \cellsize height 0pt}%
\else
\unitlength=\cellsize
\begin{picture}(1,1)
\put(0,0){\makebox(1,1){$#1$}}
\put(0,0){\line(1,0){1}}
\put(0,1){\line(1,0){1}}
\put(0,0){\line(0,1){1}}
\put(1,0){\line(0,1){1}}
\end{picture}%
\fi}}

\numberwithin{equation}{section}

\theoremstyle{plain}
\newtheorem{theorem}{Theorem}[section]
\newtheorem{proposition}[theorem]{Proposition}
\newtheorem{thmrule}[theorem]{Rule}
\newtheorem{lemma}[theorem]{Lemma}
\newtheorem{corollary}[theorem]{Corollary}
\newtheorem{conjecture}[theorem]{Conjecture}

\theoremstyle{definition}

\newtheorem{example}[theorem]{Example}

\theoremstyle{remark}
\newtheorem{remark}[theorem]{Remark}
\newtheorem{remarks}[theorem]{Remarks}

\def\fs{\mathfrak{S}}
\def\jdt{{\rm jdt}}

\begin{document}
\bibliographystyle{plain}

\begin{abstract} We present a partial  generalization to Schubert calculus on flag varieties of the classical Littlewood-Richardson rule, in its version based on Sch\"utzenberger's jeu de taquin. More precisely, we describe certain structure constants expressing the product of a Schubert and a Schur polynomial. We use a generalization of Fomin's growth diagrams (for chains in Young's lattice of partitions) to chains of permutations in the so-called $k$-Bruhat order. Our work is based on the recent thesis of Beligan, in which he generalizes the classical plactic structure on words to chains in certain intervals in $k$-Bruhat order. Potential applications of our work include the generalization of the $S_3$-symmetric Littlewood-Richardson rule due to Thomas and Yong, which is based on Fomin's growth diagrams. 
\end{abstract}

\maketitle

\section{Introduction}

Classical Schubert calculus is concerned with certain enumerative problems in geometry, which can be reduced to calculations in the cohomology of spaces such as the Grassmannian. These classical problems have been generalized in several directions, one of them being to replace the Grassmannian by the variety $Fl_n$ of   complete flags $(0=V_0\subset V_1\subset\ldots \subset V_{n-1}\subset  V_n={\mathbb C}^n)$ in ${\mathbb C}^n$. A natural basis for the cohomology of the flag variety is formed by the Schubert classes, which correspond to Schubert varieties and are indexed by permutations in the symmetric group $S_n$. The product of two Schubert classes is a positive sum of Schubert classes, as the corresponding coefficients $c_{uv}^w$ (indexed by three permutations, and known as {Schubert structure constants}) count points in a suitable triple intersection of Schubert
  varieties.

A famous open problem in algebraic combinatorics, known as the Schubert problem (and listed as Problem 11 in Stanley's survey \cite{stappc}), is to find a combinatorial description of the Schubert structure constants (and,
  in particular, a proof of their nonnegativity which bypasses geometry). The importance of this problem stems
  from the geometric significance of the Schubert structure constants, and from the fact that a
  combinatorial interpretation for these coefficients would facilitate a deeper study of
  their properties (such as their symmetries, vanishing etc.). The
  Schubert problem proved to be a very hard problem,
  resisting many attempts to be solved. The classical special case is that of the Grassmannian, in which we have the Littlewood-Richardson rule for multiplying Schur polynomials (see, e.g., \cite{staec2}[Appendix 1] or \cite{fulyt}). The corresponding coefficients $c_{\lambda\mu}^\nu$, which are indexed by three partitions, are known as Littlewood-Richardson coefficients. Currently, there are many combinatorial descriptions of these coefficients, including a recent description that reveals their $S_3$-symmetry \cite{tayssl} and is based on Fomin's growth diagrams \cite{staec2}[Appendix 1].

One usually attacks the Schubert problem for the flag variety via the multiplication of Schubert polynomials, which are polynomial representatives for Schubert classes defined by Lascoux and Sch\"utzenberger \cite{lasps}.  A notable known special case of the Schubert problem is the Pieri rule, which expresses the product of an arbitrary Schubert polynomial with one indexed by the cycle $(k+p,k+p-1,\ldots,k+1,k)$. The Pieri formula underlies the close connection between the Schubert problem and the combinatorics of chains in the so-called $k$-Bruhat order on $S_n$. We note that chains in Bruhat order are crucial objects in this area, as they underlie many of the known multiplication rules related to flag varieties (beyond the Grassmannian), including the very general formula for the $K$-theory of flag varieties of arbitrary Lie type \cite{lapawg}. We will make use of chains in $k$-Bruhat order in this paper too.

Several attempts were made to generalize the Pieri formula. Some of these attempts were based on: (1) iterating known multiplication formulas \cite{basssf,fakqad,kohmsp}; (2) proving various identities for $c_{uv}^w$ \cite{basspb,knuscr,knudcs,purvnv}; (3) bijective proofs based on insertion procedures  \cite{babrcg,koggsi,kakppf}; (4) geometric approaches \cite{bavisv,coslrr,duamrs,vakglr}. Some of the mentioned attempts led to formulas for $c_{uv}^w$ involving both positive and negative terms, but no manifestly positive formula exists in general. The most general positive rule for the Littlewood-Richardson problem is Coskun's multiplication rule for two-step flag varieties \cite{coslrr}, which is based on a geometric degeneration technique. As far as generalizing this idea to the complete flag variety $Fl_n$ is concerned, the complexity of the combinatorics involved suggests that more powerful combinatorial tools are needed.

In the joint paper with Sottile \cite{lasssp}, we defined skew Schubert polynomials based on chains in $k$-Bruhat order. The coefficients in their expansion in terms of Schubert polynomials are precisely the Schubert structure constants. Thus, we suggested an approach to the Schubert problem based on generalizing in the context of chains in Bruhat order a version of the classical Littlewood-Richardson rule which uses Sch\"utzenberger's jeu de taquin on tableaux \cite{fulyt,staec2}. The aim of the present paper is to give more details about this idea. 

A crucial piece of information underlying this paper is Beligan's recent thesis \cite{belitt} on generalizing the plactic structure for words (see, e.g., \cite{fulyt} or \cite{lltpm}) and  chains in the weak Bruhat order \cite{eagbt} to chains in $k$-Bruhat order. Beligan's results apply to maximal chains in certain intervals $[u,w]_k$ in $k$-Bruhat order which are said to contain no nesting. For such intervals, Beligan shows that each Knuth-type equivalence class of chains has as a distinguished representative a strict tableau of transpositions; this is a filling of a Young diagram with pairs $(a,b)$, such that the first entries in the pairs make the rows and columns strictly increasing. Beligan also describes an analog of the Schensted insertion algorithm (e.g., see \cite[Chapter 7.11]{staec2} or \cite{fulyt}), which transforms a chain in $[u,w]_k$ into the tableau of transpositions equivalent to it. Finally, based on this combinatorics and the Pieri formula, he generalizes this formula by describing the corresponding Schubert structure constants $c_{u,v(\lambda,k)}^w$ as the number of strict tableaux of transpositions of shape $\lambda$; here $v(\lambda,k)$ is a Grassmannian permutation corresponding to the partition $\lambda$, so the mentioned Schubert structure constants correspond to multiplying a Schubert polynomial by a Schur polynomial. In this way, Beligan generalizes the results in \cite{koggsi,kakppf}. 

In this paper, we introduce a generalization of Fomin's growth diagrams for chains in Young's lattice of partitions (which realize Sch\"utzenberger's jeu de taquin) to  chains of permutations in $k$-Bruhat order. Thus, we are able to extend the version of the classical Littlewood-Richardson rule based on jeu de taquin, which was mentioned above, to certain structure constants $c_{u,v(\lambda,k)}^w$; more precisely, we require that $u$ has no descents before or after position $k$. If we concentrate on these structure constants, which are among the ones considered by Beligan, our rule  can be viewed as a generalization of Beligan's. The reason for restricting to the above structure constants is a technical one, related to the non-nesting restriction in Beligan's work. The special case studied here gives an indication about the possible general form of our rule, which is conjectured and is currently being investigated. 

Our work has applications to the approach in \cite{lasssp} for the Schubert problem, based on skew Schubert polynomials. More importantly, it might lead to an $S_3$-symmetric description of the Schubert structure constants which would generalize the one in \cite{tayssl}. Possible connections with the geometric approaches in \cite{coslrr,vakglr} are also investigated.

\section{Background}

\subsection{The classical Littlewood-Richardson rule}\label{classlr}

Given partitions $\lambda,\mu,\nu$, let $c_{\lambda\mu}^\nu$ be the classical Littlewood-Richardson coefficient, defined as a structure constant for the multiplication of Schur functions:
\[s_\lambda\cdot s_\mu=\sum_\nu c_{\lambda\mu}^\nu s_\nu\,.\]
We briefly review the Littlewood-Richardson rule, which is a combinatorial description of $c_{\lambda\mu}^\nu$; for more details, we refer the reader to \cite{staec2}[Appendix 1] or \cite{fulyt}. 

Consider a SYT $T$ of skew shape $\nu/\mu$, and a box $b$ that can be added to $\nu/\mu$ such that the resulting shape is also a valid skew one; in addition, assume that $b$ shares its lower or right edge with $\nu/\mu$. We denote by $\jdt_b(T)$ the SYT obtained from $T$ via Sch\"utzenberger's jeu de taquin into $b$. This is given by the following simple algorithm: we pick the minimum of the entries in the boxes immediately to the right and below $b$ (there might be only one such entry) and move it to $b$; then we continue this procedure with the vacated box instead of $b$, and so on, until there is no entry to the right or below the vacated box. By applying successive jeu de taquin moves, the SYT $T$ can be transformed into a straight-shape SYT. It is well-known that the resulting SYT does not depend on the particular squence of jeu de taquin moves, so it makes sense to denote it by $\jdt(T)$. 

\begin{theorem}\label{class-lrr}(cf., e.g., {\rm \cite{staec2}[Appendix 1]}) The Littlewood-Richardson coefficient $c_{\lambda\mu}^\nu$ is equal to the number of SYT $T$ of shape $\nu/\mu$ for which $\jdt(T)$ is a particular (arbitrary) SYT $P$ of shape $\lambda$. 
\end{theorem}

We will consider two choices of $P$, which lead to two remarkable special cases of the above theorem. First, let $P_1$ be the SYT obtained by placing the entries $1,\ldots,|\lambda|$ into the boxes of $\lambda$ row by row, beginning with the top row.  Then consider all $T$ with  $\jdt(T)=P_1$. For each such $T$, define a filling $T'$ by replacing each entry $i$ in $T$ with the row number of the entry $i$ in $P$. It is not hard to see that $T'$ is a SSYT. Moreover, it is well-known that the collection of SSYT $T'$ is precisely the collection of Littlewood-Richardson tableaux of shape $\nu/\mu$ and content $\lambda$. Such a tableau is defined by the condition that its reverse row word is a lattice permutation. One also considers the so-called companion tableau of a Littlewood-Richardson tableau $T'$, which is the SSYT of shape $\lambda$ and content $\nu-\mu$ obtained by placing an entry $j$ in row $i$ of the shape $\lambda$, for each entry $i$ in row $j$ of $T'$. The companion tableaux also have a nice characterization, and a vast generalization of them is the centerpiece of Littelmann's Littlewood-Richardson rule for tensor products of irreducible representations of symmetrizable Kac-Moody algebras \cite{litlrr}. 

Now let us consider another choice $P_2$ for the SYT $P$ in Theorem \ref{class-lrr}, which is generalized by Beligan's rule for the multiplication of Schubert polynomials. We define $P_2$ as the SYT obtained by placing the entries $|\lambda|,|\lambda|-1,\ldots,1$ into the boxes of $\lambda$ column by column from right to left, going back to the rightmost unfilled column each time the first column is reached (the columns are filled from bottom to top). We will now characterize the SYT $T$ of shape $\nu/\mu$ satisfying $\jdt(T)=P_2$. For this purpose, given a SYT $T$ of shape $\nu/\mu$, we define its {\em content word} $j_1\ldots j_{|\lambda|}$ by $j_i:={\rm content}(b_i)+k$, where $b_i$ is the box containing $i$ in $T$ and $k:=\nu_1'$ (recall that the content of a box is the difference between its column and its row). It can be shown, and it is also a special case of the results below, that a SYT $T$ of shape $\nu/\mu$ satisfies $\jdt(T)=P_2$ if and only if its content word is the row word of a row and column strict tableau $T'$ of shape $\lambda$. Thus, we have the following rule.

\begin{theorem}\label{lrr2} The coefficient $c_{\lambda\mu}^\nu$ is the number of SYT of shape $\nu/\mu$ whose content word is the row word of a row and column strict tableau of shape $\lambda$.
\end{theorem}

\begin{example} We continue the example in \cite{staec2}[Appendix 1] referring to $P=P_1$, by considering $\lambda=(4,3,1)$, $\mu=(2,1)$ and $\nu=(4,4,2,1)$. We have
\[P_2=\tableau{{1}&{3}&{4}&{8}\\{2}&{6}&{7}\\{5}}\,.\]
There are two tableaux $T$, which are shown below, together with the corresponding straight-shape SYT $T'$.
\[\tableau{&&{1}&{4}\\&{3}&{7}&{8}\\{2}&{6}\\{5}}\,,\;\;\;\tableau{&&{3}&{4}\\&{1}&{7}&{8}\\{2}&{6}\\{5}}\,,\;\;\;\;\;\;\;\;\;\;\;\;\;\;\;\;\;\tableau{{1}&{3}&{5}&{6}\\{2}&{4}&{7}\\{6}}\,,\;\;\;\tableau{{1}&{3}&{5}&{6}\\{2}&{6}&{7}\\{4}}\,.\]
\end{example}

In order to better understand the definition of the content word of a SYT $T$ with shape $\nu/\mu$, let us view $T$ as a maximal chain $\mu=\mu_0<\mu_1<\ldots<\mu_{|\lambda|}=\nu$ in Young's lattice (start with $\mu$ and add to it the boxes of $\nu/\mu$ in the order indicated by the entries of $T$). Let $n:=\nu_1+\nu_1'-1=\nu_1-1+k$. To each partition $\mu_i$ corresponds a Grassmannian permutation $v_i=v(\mu_i,k)$ in $S_n$ with unique descent at $k$ (recall that this correspondence associates to the Grassmannian permutation $w=w_1\ldots w_n$ with unique descent at $k$ the partition $(w_k-k,\ldots,w_1-1)$). It turns out that we obtain a chain in the left weak Bruhat order of Grassmannian permutations $v_0<v_1<\ldots<v_{|\lambda|}$, where $v_i=s_{j_i}v_{i-1}$, with $s_{j_i}$ being the adjacent transposition $(j_i,j_i+1)$ and $j_i$ being defined as above.

We conclude this section by recalling Fomin's realization of jeu de taquin via growth diagrams. Let $T$ be a SYT of shape $\nu/\mu$, and $S$ a SYT of shape $\mu$; note that the latter tableaux determine the sequences of jeu de taquin moves for the former ones. Now consider a matrix of partitions $\lambda^{i,j}$ with $i=0,\ldots,p:=|\lambda|$ and $j=0,\ldots,q:=|\mu|$, such that its left column $(\lambda^{0,j})_j$ and its top row $(\lambda^{i,q})_i$ are the maximal chains in Young's lattice corresponding to $S$ and $T$, respectively; in particular, $\lambda^{0,0}$ is the empty partition, $\lambda^{0,q}=\mu$, and $\lambda^{p,q}=\nu$. The other partitions are determined by the local rule below, which specifies $\lambda^{i+1,j}$ based on $\lambda^{i,j}$, $\lambda^{i,j+1}$, and $\lambda^{i+1,j+1}$. In order to specify this rule, let 
\begin{equation}\label{defhv}\sigma_{i,j}:=\lambda^{i+1,j}-\lambda^{i,j}\;\;\;\;\mbox{and}\;\;\;\;\tau_{i,j}:=\lambda^{i,j+1}-\lambda^{i,j}\,.\end{equation}
We call $\sigma_{i,j}$ and $\tau_{i,j}$ horizontal and vertical transpositions, respectively.

\begin{thmrule}\label{fominrule} If the two boxes of $\lambda^{i+1,j+1}\setminus \lambda^{i,j}$ are not adjacent, then $\sigma_{i,j}=\sigma_{i,j+1}$ ${\rm (}$and $\tau_{i+1,j}=\tau_{i,j}{\rm )}$; otherwise, $\sigma_{i,j}=\tau_{i,j}$ ${\rm (}$and $\tau_{i+1,j}=\sigma_{i,j+1}{\rm )}$. 
\end{thmrule}

Note that in the first case of the rule, the interval $[\lambda^{i,j},\lambda^{i+1,j+1}]$ in Young's lattice is a product of two chains of length 1, whereas in the second case it is a chain of length 2. It turns out that the bottom row of the above matrix of partitions, namely $(\lambda^{i,0})_i$, is precisely the chain corresponding to $\jdt(T)$. Thus, the SYT $T$ in Theorem \ref{class-lrr} are determined by the condition that $(\lambda^{i,0})_i$ corresponds to the fixed SYT $P$ (in particular, $\lambda^{p,0}=\lambda$). 

\subsection{Beligan's rule for multiplying Schubert polynomials}\label{beligan} Given a permutation $w$ in the symmetric group $S_n$, we denote by $\fs_w(x)$ the Schubert polynomial indexed by $w$. This is a homogeneous polynomial in $x_1,\ldots,x_{n-1}$ with nonnegative integer coefficients and degree $\ell(w)$ (the length of $w$). For more information on Schubert polynomials, we refer the reader to \cite{fulyt,lasps,macnsp,manfsp}.

The main outstanding problem in the theory of Schubert polynomials is the
Littlewood-Richardson problem~\cite[Problem 11]{stappc}:
Determine the structure constants
$c^w_{uv}$ defined by the polynomial identity
 \[
    \fs_u(x) \cdot \fs_v(x)\ =\ \sum_w c^w_{uv} \fs_w(x)\,.
 \]
Since every Schur polynomial is a Schubert polynomial, this problems asks for
the generalization of the classical Littlewood-Richardson rule.
The Littlewood-Richardson coefficients $c^w_{uv}$ for
Schubert polynomials are important since they are intersection numbers of Schubert varieties; more
precisely, $c^w_{uv}$ enumerates flags in a suitable triple
intersection of  Schubert varieties indexed by $u$, $v$, and $w_\circ w$ (where $w_\circ$ is the longest permutation in $S_n$).

Beligan \cite{belitt} gave a combinatorial interpretation for certain coefficients $c_{uv}^w$ when $u=u(\lambda,k)$ is a Grassmannian permutation (with unique descent at $k$, cf. the previous section). In other words, this rules gives certain coefficients in the multiplication of a Schubert polynomial by a Schur polynomial. It is known from \cite{basspb} that $c_{u(\lambda,k),v}^w=0$ unless $v<w$ in the so-called {\em $k$-Bruhat order} on the symmetric group $S_n$ and $\ell(w)-\ell(v)=|\lambda|$. This highlights the importance of the $k$-Bruhat order, which we now recall. The  Bruhat order is the partial order on $S_n$ with covering relations $v\lessdot w=v(a,b)$, where
$\ell(w)=\ell(v)+1$ and $(a,b)$ denotes the transposition of $a<b$. A permutation $v$ admits a cover $v\lessdot v(a,b)$ with $a<b$ and
$v(a)<v(b)$ if and only if whenever $a<c<b$, then either $v(c)<v(a)$ or else
$v(b)<v(c)$. This is known as the {\it cover condition}; it is both explicitly and implicitly used several times in this paper. The $k$-Bruhat order, denoted $<_k$, is the suborder of the Bruhat order where the covers are restricted to those $v\lessdot v(a,b)$ with $a\le k<b$.

A crucial role in Beligan's rule is played by maximal chains in $k$-Bruhat order $w_0\lessdot w_1\lessdot w_2\lessdot \ldots$. We denote these chains by words of transpositions $\alpha_\beta \gamma_\delta \ldots$, where $w_1=(\alpha,\beta)w_0$, $w_2=(\gamma,\delta)w_1$, while $\alpha<\beta$, $\gamma<\delta$ etc. Note that, as opposed to the definition of the $k$-Bruhat order, here we use left multiplication. An interval in the $k$-Bruhat order, denoted $[v,w]_k$, is said to contain  {\em no nesting} if none of its maximal chains contains a segment of the form $\alpha_\delta\beta_\gamma$ or $\beta_\gamma\alpha_\delta$, where $\alpha<\beta<\gamma<\delta$. Beligan's rule refers to those structure constants $c_{u(\lambda,k),v}^w$ for which $[v,w]_k$ contains no nesting. Some criteria for the non-nesting property of $[v,w]_k$ are given in \cite{belitt}. The only ones that involve only $v$ or only $w$ seem to be the following:
\begin{enumerate}
\item $v$ has no descents after position $k$ (i.e., it is a {\em $k$-semi-shuffle});
\item $v$ has no descents before position $k$;
\item $w$ has no ascents after position $k$;
\item $w$ has no ascents before position $k$.
\end{enumerate}
The following simple lemma about intervals containing no nesting will be useful.

\begin{lemma}\label{simplelem} If $\alpha_\beta\gamma_\delta$ is a subchain of a maximal chain in an interval in $k$-Bruhat order containing no nesting, then the following are equivalent: (i) $\alpha<\gamma$; (ii) $\beta<\delta$; (iii) $\beta\le \gamma$.
\end{lemma}

In order to state Beligan's rule, we need to consider tableaux of transpositions, that is, fillings of Young diagrams with transpositions. If the first entries in the transpositions are (strictly) increasing in rows and columns, the tableau is called {\em strict}. The row word of a tableau is defined as usual. 

\begin{theorem}\cite{belitt}\label{thmbel} If $[v,w]_k$ contains no nesting, then $c_{u(\lambda,k),v}^w$ is equal to the number of strict tableaux of transpositions of shape $\lambda$ whose row word is a maximal chain in $[v,w]_k$. 
\end{theorem}

\begin{remark} Let us assume that $v$ and $w$ are Grassmannian permutations with unique descents at $k$, that is, $v=v(\mu,k)$ and $w=w(\nu,k)$. A maximal chain in $k$-Bruhat order from $v$ to $w$ (which is, in fact, a chain in the left weak Bruhat order) corresponds to a maximal chain in Young's lattice from $\mu$ to $\nu$, that is, to a SYT of shape $\nu/\mu$. Thus, Theorem \ref{lrr2} is a special case of Theorem \ref{thmbel}.
\end{remark}

A crucial ingredient in the proof of the above theorem is the so-called plactic structure on the maximal chains in intervals $[v,w]_k$ containing no nesting (the latter condition is implicit from now on). The mentioned plactic structure generalizes the classical plactic structure on words (see, e.g., \cite{fulyt}), which is relevant to the classical Littlewood-Richardson rule. The plactic relations for maximal chains in $k$-Bruhat order are the following:
\begin{enumerate}
\item[(KB1)] $\;\;\;\alpha_\gamma\gamma_\delta\beta_\gamma\sim \beta_\gamma\alpha_\beta\beta_\delta\,,$
\item[(KB2)] $\;\;\;\beta_\gamma \gamma_\delta \alpha_\gamma \sim \beta_\delta\alpha_\beta\beta_\gamma\,,\;\;\;\;\;\;$ where $\alpha<\beta<\gamma<\delta$;
\item[(KB3)] $\;\;\;\alpha_\beta\varepsilon_\varphi\gamma_\delta\sim \varepsilon_\varphi\alpha_\beta\gamma_\delta\,,$
\item[(KB4)] $\;\;\;\gamma_\delta\varepsilon_\varphi\alpha_\beta\sim \gamma_\delta\alpha_\beta\varepsilon_\varphi\,,\;\;\;\;\;\;$ where $\alpha<\beta\le\gamma<\delta\le\varepsilon<\varphi$.
\end{enumerate}

As usual, one defines the plactic equivalence of maximal chains in $[v,w]_k$ as the equivalence relation generated by (KB1)-(KB4). Then we have the following generalization of the corresponding classical result.

\begin{theorem}\cite{belitt}\label{plac1}
Each plactic equivalence class contains a unique chain which is the row word of a strict tableau of transpositions. 
\end{theorem} 

As in the classical case, the tableau $P(\Gamma)$ equivalent to a chain $\Gamma$ can be obtained from $\Gamma$ by an insertion procedure, which is described in detail in \cite{belitt}. It is shown that this procedure can be reduced to a sequence of substitutions of the right-hand sides of (KB1)-(KB4) for the corresponding left-hand sides. Moreover, we can form a recording SYT $Q(\Gamma)$ of the same shape as $P(\Gamma)$, as usual. As expected, we have the following generalization of the corresponding classical result.

\begin{theorem}\cite{belitt}\label{plac2}
The correspondence $\Gamma\mapsto (P(\Gamma),Q(\Gamma))$ is a bijection between maximal chains in $[v,w]_k$ and pairs consisting of a strict tableau of transpositions (for a maximal chain in $[v,w]_k$) and a SYT of the same shape.
\end{theorem}

\section{Main results and conjectures}

We now define our new jeu de taquin on maximal chains in $k$-Bruhat order, which generalizes Sch\"utzenberger's jeu de taquin in Fomin's realization, based on growth diagrams.

We begin with an arbitrary maximal chain $\Delta=(v_0\lessdot_l v_1\lessdot_l\ldots\lessdot_l v_q=v)$ in $l$-Bruhat order starting at the identity, followed by a maximal chain $\Gamma=(v=w_0\lessdot_k w_1\lessdot_k\ldots\lessdot_k w_p=w)$ in $k$-Bruhat order. As in the classical case, we form a matrix of permutations $w^{i,j}$ given by a local rule to be specified and the following boundary conditions:
\[w^{0,j}:=v_j\,,\;\;\;\;\;\; w^{i,q}:=w_i\,,\;\;\;\;\;\;\;\;\;\mbox{for $i=0,\ldots,p\;\;$ and $\;j=0,\ldots,q$}\,.\]
We then set
\[\jdt_\Delta(\Gamma):=(w^{0,0},w^{1,0},\ldots,w^{p,0})\,.\]

Let us now specify the local rule, which amounts to specifying $w^{i+1,j}$ based on $w^{i,j}$, $w^{i,j+1}$, and $w^{i+1,j+1}$. Note that the interval $[w^{i,j},w^{i+1,j+1}]$ in Bruhat order is always a product of two chains of length 1. Let us denote  by $x$ the unique permutation in the corresponding open interval that is different from $w^{i,j+1}$. It is possible define the local rule as follows.

\begin{thmrule}\label{genrule} Set $w^{i+1,j}:=x$ or $w^{i+1,j}:=w^{i,j+1}$ such that $w^{i,j}\lessdot_k w^{i+1,j}$ and $w^{i+1,j}\lessdot_l w^{i+1,j+1}$. If both choices work, give preference to the first one. 
\end{thmrule}

Thus, all the horizontal chains are in $k$-Bruhat order, while all the vertical ones are in $l$-Bruhat order; in particular, $\jdt_\Delta(\Gamma)$ is a chain in $k$-Bruhat order. Let us now discuss in detail all the cases that can appear, in order to justify the claim that the rule can always be applied; we will see that, compared to the classical (Grassmannian) case, we have an increase from 2 (in fact, 3) to 13 cases. Recalling the notation for chains as words of transpositions which are applied on the left, we define the transpositions $\sigma_{i,j}$ and $\tau_{i,j}$ as in (\ref{defhv}):
\begin{equation}\sigma_{i,j}:=w^{i+1,j}(w^{i,j})^{-1}\;\;\;\;\mbox{and}\;\;\;\;\tau_{i,j}:=w^{i,j+1}(w^{i,j})^{-1}\,.\end{equation}
We will specify the 13 cases as $\tau_{i,j}\sigma_{i,j+1}\mapsto \sigma_{i,j}\tau_{i+1,j}$ (cf. Rule \ref{fominrule}), depending on some condition on $w=w^{i,j}$. In all the cases except the first one, it is assumed that $\alpha<\beta<\gamma$.

\begin{enumerate}
\item[(J0)] $\;\;\;\alpha_\beta\gamma_\delta\mapsto\gamma_\delta\alpha_\beta$ if $\alpha,\beta,\gamma,\delta$ are all distinct;
\item[(J1)] $\;\;\;\alpha_\beta\alpha_\gamma\mapsto\beta_\gamma\alpha_\beta$;
\item[(J2)] $\;\;\;\alpha_\gamma\alpha_\beta\mapsto\alpha_\beta\beta_\gamma$;
\item[(J3)] $\;\;\;\beta_\gamma\alpha_\gamma\mapsto\alpha_\beta\beta_\gamma$;
\item[(J4)] $\;\;\;\alpha_\gamma\beta_\gamma\mapsto\beta_\gamma\alpha_\beta$;
\item[(J5)] $\;\;\;\alpha_\beta\beta_\gamma\mapsto\beta_\gamma\alpha_\gamma$ if $w^{-1}(\beta)\le k<w^{-1}(\gamma)$;
\item[(J5$'$)] $\;\;\;\alpha_\beta\beta_\gamma\mapsto\alpha_\beta\beta_\gamma$ if $k<w^{-1}(\beta)<w^{-1}(\gamma)$;
\item[(J6)] $\;\;\;\alpha_\beta\beta_\gamma\mapsto\alpha_\gamma\alpha_\beta$ if $w^{-1}(\gamma)\le l<w^{-1}(\beta)$;
\item[(J6$'$)] $\;\;\;\alpha_\beta\beta_\gamma\mapsto\alpha_\beta\beta_\gamma$ if $l<w^{-1}(\gamma)<w^{-1}(\beta)$;
\item[(J7)] $\;\;\;\beta_\gamma\alpha_\beta\mapsto\alpha_\beta\alpha_\gamma$ if $w^{-1}(\alpha)\le k<w^{-1}(\beta)$;
\item[(J7$'$)] $\;\;\;\beta_\gamma\alpha_\beta\mapsto\beta_\gamma\alpha_\beta$ if $w^{-1}(\alpha)<w^{-1}(\beta)\le k$;
\item[(J8)] $\;\;\;\beta_\gamma\alpha_\beta\mapsto\alpha_\gamma\beta_\gamma$ if $w^{-1}(\beta)\le l<w^{-1}(\alpha)$;
\item[(J8$'$)] $\;\;\;\beta_\gamma\alpha_\beta\mapsto\beta_\gamma\alpha_\beta$ if $w^{-1}(\beta)<w^{-1}(\alpha)\le l$.
\end{enumerate}

\begin{remarks}\label{3rem} (1) By inspecting the above cases, we can see that there is a unique choice in Rule \ref{genrule} as long as $\tau_{i,j}$ and $\sigma_{i,j+1}$ do not commute. Otherwise, we can have one or two choices.

(2) The Grassmannian cases (when $k=l$) fall under (J0), (J5$'$) $-$ the two adjacent boxes are in the same row, and (J7$'$) $-$ the two adjacent boxes are in the same column (cf. Rule \ref{fominrule}). 

(3) The above jeu de taquin can be generalized by letting $\Delta$ and $\Gamma$ be concatenations of maximal chains in $k$-Bruhat order for various $k$; such chains will be called mixed Bruhat chains. This more general version of jeu de taquin will be needed below.
\end{remarks}

\begin{example}\label{grdiag}
Let $w^{0,0}:=2143$, $l=1$, $k=2$, $\Delta=2_4$, $\Gamma=1_22_3$. We have $\jdt_\Delta(\Gamma)=1_42_3$, where the applied transformations are (J8) and (J2). We can indicate the growth diagram pictorially, as follows (the arrows indicate the increasing direction in Bruhat order).
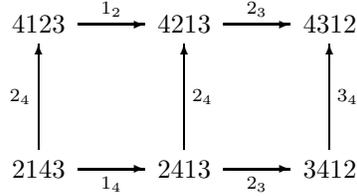
\begin{equation*} 
\begin{diagram}
\node{4123}
\arrow{e,t}{1_2}\node{4213}\arrow{e,t}{2_3}\node{4312}\\
\node{2143}\arrow{n,l}{2_4}\arrow{e,b}{1_4}\node{2413}\arrow{n,r}{2_4}\arrow{e,b}{2_3}\node{3412}\arrow{n,r}{3_4}
\end{diagram}
\end{equation*}
\end{example}

By Remark \ref{3rem} (1), our jeu de taquin for chains in $k$-Bruhat order is symmetric, in the sense stated below.

\begin{proposition} Consider the input of jeu de taquin to be the pair $(\Delta,\Gamma)$ and the output the pair $(\Gamma',\Delta')$, where $\Gamma':=\jdt_\Delta(\Gamma)$ and $\Delta'=(w^{p,0}\lessdot_l w^{p,1}\lessdot_l\ldots\lessdot_l w^{p,q})$. Then  if we input $(\Gamma',\Delta')$, the output is $(\Delta,\Gamma)$. \end{proposition}

We now state our main conjecture, which is the natural generalization of Theorem \ref{class-lrr}, by Remark \ref{3rem} (2).

\begin{conjecture}\label{lrr} Consider permutations $v\le_k w$ and a Grassmannian permutation $u(\lambda,k)$. There is a mixed Bruhat chain $\Delta$ from the identity to $v$ such that the Littlewood-Richardson coefficient $c_{u(\lambda,k),v}^w$ is equal to the number of maximal chains $\Gamma$ in $k$-Bruhat order from $v$ to $w$ for which $\jdt_\Delta(\Gamma)$ is a particular (arbitrary) maximal chain $\Gamma'$ in $k$-Bruhat order from the identity to $u(\lambda,k)$. 
\end{conjecture}

\begin{remarks}\label{remsyt} (1) In Conjecture \ref{lrr}, $\Gamma'$ can be thought of as a maximal chain in Young's lattice from the empty partition to $\lambda$, or simply as a SYT of shape $\lambda$. 

(2) It is not true that all mixed Bruhat chains from the identity to $v$ satisfy the condition in Conjecture \ref{lrr}. The reasons for this will be discussed in detail in Section \ref{prooflem}. This situation is different from the Grassmannian case, where the jeu de taquin moves can be performed in any order (i.e., along any chain in Young's lattice). 
\end{remarks}

We will now prove two special cases of this conjecture, which should give an idea about the general case. Our work is based on the results of Beligan about the plactic structure of intervals in $k$-Bruhat order containing no nesting. The main problem in attacking the general case is that the plactic structure is only understood in the non-nesting case so far, while the fact that $[v,w]_k$ contains no nesting does not guarantee that $\jdt_\Delta(\Gamma)$ is a maximal chain in an interval with the same property (say, for a chain $\Delta$ of length 1, as Example \ref{grdiag} shows). On the other hand, the only known simple criteria for the non-nesting property are those in Section \ref{beligan}. Therefore, our main results are related to the cases when $v$ has no descents before or after position $k$, and are based on mixed Bruhat chains $\Delta=(v_0\lessdot v_1\lessdot\ldots\lessdot v_q=v)$ with the following properties. 
\begin{enumerate}
\item[(PL)] For each $i$, $v_i$ is a $k$-semi-shuffle. Furthermore, if $v_{i+1}=v_i(l,\,\cdot\,)$, then $l$ is the smallest non-fixed point of $v_{i+1}$, and $v_i\lessdot_l v_{i+1}$. 
\item[(PR)] For each $i$, $v_i$ has no descents before $k$. Furthermore, if $v_{i+1}=v_i(\,\cdot\,,l)$, then $l$ is the largest non-fixed point of $v_{i+1}$, and $v_i\lessdot_{l-1} v_{i+1}$. 
\end{enumerate}

\begin{remarks}\label{remplr} (1) Clearly, if $v$ is a $k$-semi-shuffle, then there exists a mixed Bruhat chain from the identity to $v$ with property (PL). Moreover, this chain is a concatenation $\Delta^{k}\ldots\Delta^{1}$, where $\Delta^l$ is a chain in $l$-Bruhat order. Similarly, if $v$ has no descents before position $k$, then there exists a mixed Bruhat chain from the identity to $v$ with property (PR), and this chain is a concatenation $\Delta^{k+1}\ldots\Delta^{n-1}$.

(2) If the chain $\Delta$ in a growth diagram satisfies property (PL), resp. (PR), then all permutations in the growth diagram are $k$-semi-shuffles, resp. have no descents before $k$. This remark will be used implicitly below.
\end{remarks}

With this notation, we can state our main result.

\begin{theorem}\label{speclrr} {\rm (1)} If $v$ is a $k$-semi-shuffle, then any mixed Bruhat chain $\Delta$ from the identity to $v$ with property {\rm (PL)} satisfies the condition in Conjecture {\rm \ref{lrr}}.

{\rm (2)} The same is true if $v$ has no descents before $k$ for $\Delta$ having property {\rm (PR)}.

In both cases, $\jdt_\Delta(\Gamma)$ does not depend on the chains $\Delta$ having the mentioned properties. 
\end{theorem}

Remark \ref{remplr} (1) leads us to defining property (PLR) as the natural generalization of properties (PL) and (PR) for chains $\Delta$. We use the above notation. 
\begin{enumerate}
\item[(PLR)] The chain $\Delta$ is a concatenation $\Delta'\Delta''$ (resp. $\Delta''\Delta'$), where $\Delta'$ has property (PL) and $\Delta''=\Delta^{k+1}\ldots\Delta^{n-1}$ (resp. $\Delta''$ has property (PR) and $\Delta'=\Delta^{k}\ldots\Delta^{1}$). Furthermore, given $v_i\lessdot_l v_{i+1}$ in $\Delta''$ (resp. $\Delta'$), we have $v_{i+1}=v_i(\,\cdot\,,l+1)$ (resp. $v_{i+1}=v_i(l,\,\cdot\,)$). 
\end{enumerate} 

\begin{remark} For any permutation $v$, there exists a mixed Bruhat chain from the identity to $v$ with property (PLR). Indeed, the goal is to go down in Bruhat order from $v$ along the reverse of a chain $\Delta''$ (resp. $\Delta'$) with the properties above, until we reach a $k$-semi-shuffle (resp. a permutation with no descents before position $k$). 
\end{remark}

Theorem \ref{speclrr} suggests the following stronger version of Conjecture \ref{lrr}, which we currently investigate, and which was so far confirmed by several computer tests. 

\begin{conjecture} There exists a mixed Bruhat chain $\Delta$ from the identity to $v$ with property {\rm (PLR)} which satisfies the condition in Conjecture {\rm \ref{lrr}}.
\end{conjecture}

Theorem \ref{speclrr} (and, in fact, the Conjecture \ref{lrr} in general) can be proved based on the following result.

\begin{proposition}\label{presq} A chain $\Delta$ satisfies the condition in Conjecture {\rm \ref{lrr}} if 
\begin{equation}\label{eqpresq}Q(\jdt_\Delta(\Gamma))=Q(\Gamma)\end{equation}
for any maximal chain $\Gamma$ in $[v,w]_k$. In this case, $\jdt_\Delta(\Gamma)$ does not depend on $\Delta$. 
\end{proposition}

\begin{proof} Clearly, a maximal chain in $k$-Bruhat order starting at the identity consists only of Grassmannian permutations. Moreover, it is well-known that all maximal chains in an interval $[1,u(\lambda,k)]_k$ have the same $P$-tableau, which is of shape $\lambda$. (In fact only the relations (KB3) and (KB4) are needed here to get the plactic equivalence.) So a chain in the mentioned interval is determined by its $Q$-tableau. If an arbitrary maximal chain $\Gamma'$ in $[1,u(\lambda,k)]_k$ is fixed, then, by (\ref{eqpresq}) and Theorems \ref{plac1}-\ref{plac2}, the number of maximal chains $\Gamma$ in $[v,w]_k$ satisfying $\jdt_\Delta(\Gamma)=\Gamma'$ is just the number of plactic classes in $[v,w]_k$ represented by tableaux of shape $\lambda$. But the latter number is $c_{u(\lambda,k),v}^w$, by Theorem \ref{thmbel}. 
\end{proof}

\begin{remark} The proof of Proposition \ref{presq} clarifies the way in which Theorem \ref{speclrr} is a generalization of Theorem \ref{thmbel} for the considered permutations $v$. More precisely, the maximal chains in the latter theorem are are precisely the chains $\Gamma$ in Conjecture \ref{lrr} when $\Gamma'$ is the SYT $P_2$ considered in Section \ref{classlr} (upon the identification of chains and SYT in Remark \ref{remsyt} (1)).
\end{remark}

The next lemma is our main tool for proving Theorem \ref{speclrr} via Proposition \ref{presq}. The lemma will be proved in the following section based on a detailed case by case analysis of the interaction between the plactic relations (KB1)-(KB4) and the local rules (J0)-(J8$'$) for growth diagrams.

\begin{lemma}\label{presplac}
Let $(\overline{\Delta},\overline{\Gamma})\mapsto (\overline{\Gamma}',\overline{\Delta}')$ be a height $1$ fragment of a growth diagram for which $\Delta$  has properties {\rm (PL)} or {\rm (PR)}. If $\overline{\Gamma}$ is the left-hand side of a plactic relation, then so is $\Gamma'$; moreover, if $\overline{\Gamma}\sim\overline{\Gamma}_r$ and $\overline{\Gamma}'\sim\overline{\Gamma}_r'$, then $(\overline{\Delta},\overline{\Gamma}_r)\mapsto(\overline{\Gamma}_r',\overline{\Delta}')$. On the other hand, if $\overline{\Gamma}=\alpha_\beta\gamma_\delta$ and $\overline{\Gamma}'=\alpha'_{\beta'}\gamma'_{\delta'}$, then $\alpha<\gamma$ if and only if $\alpha'<\gamma'$. 
\end{lemma}

We have the following corollary to the above lemma.

\begin{corollary}\label{corpres} Let $({\Delta},{\Gamma})\mapsto ({\Gamma}',{\Delta}')$ be a growth diagram for which $\Delta$  has properties {\rm (PL)} or {\rm (PR)}. Then we have
\[\jdt_\Delta(P(\Gamma))=P(\Gamma')\,,\;\;\;\;\mbox{and}\;\;\;\;Q(\Gamma)=Q(\Gamma')\,.\]
\end{corollary}

\begin{proof} It suffices to consider $\Delta$ of length 1. Recall that the insertion algorithm in \cite{belitt} can be reduced to successively replacing the left-hand sides of relations (KB1)-(KB4) with the corresponding right-hand sides; more precisely, a nontrivial insertion in a row of length $n$ amounts to applying the mentioned relations in positions $n-1,n-2,\ldots,1$, just like in the usual Schensted insertion $-$ see \cite{belitt}[Section 4.3]. Based on this and Lemma \ref{presplac}, the corollary follows. Indeed, Lemma \ref{presplac} says that each time we apply a plactic relation in the insertion algorithm for the top row of a growth diagram, we can apply the corresponding plactic relation in the bottom row, and still have a growth diagram. The chains $\overline{\Gamma}$ of length 2 in Lemma \ref{presplac} are needed to take care of the usual comparison test for inserting a letter at the end of a row.  
\end{proof}

\begin{proof}[Proof of Theorem {\rm \ref{speclrr}}]  Immediate by Corollary \ref{corpres} and Proposition \ref{presq}. 
\end{proof}

\section{Proof of Lemma \ref{presplac}}\label{prooflem}

The goal is to verify all the instances of Lemma \ref{presplac} in the two cases considered, namely when $\Delta$ has property (PL) and (PR). We start with some simple lemmas. 

\begin{lemma}\label{ruleoutkb} Relation {\rm (KB1)} can be applied only to a permutation for which $\beta$ precedes $\alpha$, both entries being in positions $1$ through $k$. Relation {\rm (KB2)} can be applied only to a permutation for which $\delta$ precedes $\gamma$, both entries being in positions larger than $k$.
\end{lemma}

\begin{proof} This is immediate by the cover condition.
\end{proof}

We first concentrate on property (PL). By Lemma \ref{ruleoutkb}, only the plactic relations (KB1), (KB3), and (KB4) can be applied in this case. 

\begin{lemma}\label{ruleout} If a chain $\Delta$ has property {\rm (PL)}, then rules {\rm (J3), (J4), (J6), (J6$'$), (J7)}, and {\rm (J8$'$)} are never applied.
\end{lemma}

\begin{proof}
For the first five rules, the statement is clear because some permutation in the matrix of permutations would not be a $k$-semi-shuffle. Assume that (J8$'$) is applied and that $\Delta$ has length 1; thus $\tau_{i,0}\sigma_{i,1}=\sigma_{i,0}\tau_{i+1,0}=\beta_\gamma\alpha_\beta$ with $\alpha<\beta<\gamma$, for some $i$. Let $r:=(w^{i,0})^{-1}(\beta)$ and $s:=(w^{i,0})^{-1}(\alpha)$, where $r<s$, by the condition for rule (J8$'$). The vertical chains are in $l$-Bruhat order, where $l\ge s$, again by the condition for rule (J8$'$); in fact, $r<s\le l\le k$. We have $w^{0,0}(r)>w^{0,0}(s)$, since inversions of entries in positions 1 through $k$ are preserved upon going down in a maximal chain in $k$-Bruhat order.  Since $\tau_{0,0}$ exchanges the entry in position $l$ of $w^{0,0}$ with an entry to its right, while a previous entry (in position $r$) is inverted with an entry to its right, property (PL) is contradicted.
\end{proof}

The next lemma provides the main criterion for ruling out the cases that contradict Lemma \ref{presplac}.

\begin{lemma}\label{movevert} Consider $\Delta$ having property {\rm (PL)}, and assume it has length $1$. Then it is not possible to have $\tau_{i,0}=(\alpha,\,\cdot\,)$ for some $i$, and an entry $\beta>\alpha$ in $w^{i,0}$ to the left of $\alpha$. 
\end{lemma}

\begin{proof} Let $r:=(w^{i,0})^{-1}(\beta)$ and $s:=(w^{i,0})^{-1}(\alpha)$, where $r<s$. Upon a case by case examination, based on Lemma \ref{ruleout}, we see that the left entries in the vertical transpositions $\tau_{j,0}$ with $j\le i$ are in positions greater or equal to $s$. Thus $l\ge s$. We now have an identical situation with the one in the proof of Lemma \ref{ruleout}, so property (PL) is contradicted.
\end{proof}

\subsection{Preserving monotonicity}

From now on we will use the notation in Lemma \ref{presplac} freely. We start by examining 2 by 1 fragments of  growth diagrams. The following remark will be useful throughout.

\begin{remark}\label{simplerem} If $\alpha_\beta\gamma_\delta\mapsto\gamma'_{\delta'}\alpha'_{\beta'}$ by one of the rules (J0)-(J8$'$), then each of the intersections $\{\alpha,\beta\}\cap\{\alpha',\beta'\}$ and $\{\gamma,\delta\}\cap\{\gamma',\delta'\}$ contains at least an element. This remark will be used implicitly, often in conjunction with Lemma \ref{simplelem}.
\end{remark}

\begin{lemma}\label{less} Let $\overline{\Gamma}=\alpha_\beta\gamma_\delta$ and $\overline{\Gamma}'=\alpha'_{\beta'}\gamma'_{\delta'}$. If $\alpha<\gamma$, then $\alpha'<\gamma'$.
\end{lemma}

\begin{proof}
By Remark \ref{simplerem}, it suffices to consider the case $\beta=\gamma$. Assume that $\alpha'>\gamma'$. This can only happen if $\alpha'=\delta'=\beta$. By inspecting the jeu de taquin rules, we conclude that the fragment of the growth diagram has the following form.

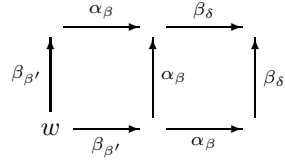
\begin{equation*} 
\begin{diagram}
\node{}
\arrow{e,t}{\alpha_\beta}\node{}\arrow{e,t}{\beta_\delta}\node{}\\
\node{w}\arrow{n,l}{\beta_{\beta'}}\arrow{e,b}{\beta_{\beta'}}\node{}\arrow{n,r}{\alpha_\beta}\arrow{e,b}{\alpha_\beta}\node{}\arrow{n,r}{\beta_\delta}
\end{diagram}
\end{equation*}

In other words, the two rules applied are (J7$'$) and (J5$'$). The entries $\alpha$ and $\beta$ are in positions 1 through $k$ of the permutation $w$ in the bottom left corner of the growth diagram, while $\beta'$ and $\delta$ are in positions larger than $k$. It is easy to see that we cannot have $\beta'=\delta$. If $\beta'>\delta$, then $\beta'$ is to the right of $\delta$ in $w$ (since $w$ is a $k$-semi-shuffle), and the transposition $(\beta,\beta')$ applied to $w$ violates the cover condition. Thus, we must have $\beta'<\delta$. By rule (J7$'$), $\alpha$ precedes $\beta$ in $w$. Since we have the chain in $k$-Bruhat order $\beta_{\beta'}\alpha_\beta\beta_\delta$, Lemma \ref{ruleoutkb} tells us that $\beta$ precedes $\alpha$ in $w$, which is a contradiction.
\end{proof}

\begin{lemma}\label{gtr} Let $\overline{\Gamma}=\alpha_\beta\gamma_\delta$ and $\overline{\Gamma}'=\alpha'_{\beta'}\gamma'_{\delta'}$. If $\alpha>\gamma$, then $\alpha'>\gamma'$.
\end{lemma}

\begin{proof} Like in the proof of the previous lemma, it suffices to consider a special case, namely $\alpha=\delta$. Assume that $\alpha'<\gamma'$. This can only happen if $\beta'=\gamma'=\alpha$. By inspecting the jeu de taquin rules, we conclude that the fragment of the growth diagram has the following form.

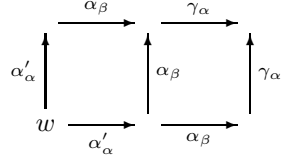
\begin{equation*} 
\begin{diagram}
\node{}
\arrow{e,t}{\alpha_\beta}\node{}\arrow{e,t}{\gamma_\alpha}\node{}\\
\node{w}\arrow{n,l}{\alpha'_{\alpha}}\arrow{e,b}{\alpha'_{\alpha}}\node{}\arrow{n,r}{\alpha_\beta}\arrow{e,b}{\alpha_\beta}\node{}\arrow{n,r}{\gamma_\alpha}
\end{diagram}
\end{equation*}

In other words, the two rules applied are (J5$'$) and (J7$'$). The entries $\alpha'$ and $\gamma$ are in positions 1 through $k$ of the permutation $w$ in the bottom left corner of the growth diagram, while $\alpha$ and $\beta$ are in positions larger than $k$. It is easy to see that we cannot have $\alpha'=\gamma$. The case $\alpha'>\gamma$ is ruled out by Lemma \ref{ruleoutkb}, since we have the chain $\alpha'_\alpha\alpha_\beta\gamma_\alpha$ in $k$-Bruhat order starting at $w$.  Thus, we must have $\alpha'<\gamma$. By rule (J7$'$), $\gamma$ precedes $\alpha'$ in $w$. Then we obtain a contradiction by Lemma \ref{movevert}.
\end{proof}

Lemmas \ref{less} and \ref{gtr} verify Lemma \ref{presplac} for 2 by 1 fragments of growth diagrams. 

\subsection{The relation (KB3)}

The following lemma prepares the case when $\overline{\Gamma}$ is the left-hand side of a relation (KB3).

\begin{lemma}\label{lesstoeq}
Let $\overline{\Gamma}=\alpha_\beta\gamma_\delta$ and $\overline{\Gamma}'=\alpha'_{\beta'}\gamma'_{\delta'}$, where $\beta<\gamma$. The only cases when $\beta'=\gamma'$ are the ones shown in the diagrams below; in the first case, $w$ has the form $\ldots\alpha\ldots\beta\ldots|\ldots\gamma\ldots\delta\ldots$, while in the second one it has the form $\ldots\beta\ldots\alpha\ldots|\ldots\gamma\ldots\delta\ldots$ (the vertical bar is between positions $k$ and $k+1$).
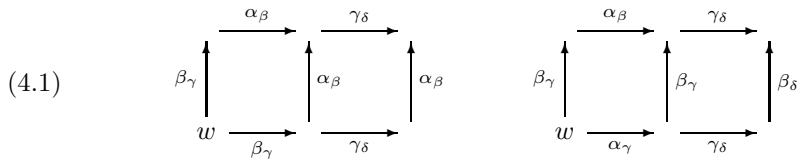
\begin{equation}\label{2grd} 
\begin{diagram}
\node{}
\arrow{e,t}{\alpha_\beta}\node{}\arrow{e,t}{\gamma_\delta}\node{}\\
\node{w}\arrow{n,l}{\beta_\gamma}\arrow{e,b}{\beta_\gamma}\node{}\arrow{n,r}{\alpha_\beta}\arrow{e,b}{\gamma_\delta}\node{}\arrow{n,r}{\alpha_\beta}
\end{diagram}\;\;\;\;\;\;\;\;\;\;\;\;\;\;\;\;
\begin{diagram}
\node{}
\arrow{e,t}{\alpha_\beta}\node{}\arrow{e,t}{\gamma_\delta}\node{}\\
\node{w}\arrow{n,l}{\beta_\gamma}\arrow{e,b}{\alpha_\gamma}\node{}\arrow{n,r}{\beta_\gamma}\arrow{e,b}{\gamma_\delta}\node{}\arrow{n,r}{\beta_\delta}
\end{diagram}
\end{equation}
If $\beta'<\gamma'$, then there is only one case when none of the two jeu de taquin rules is {\rm (J0)}, namely the growth diagram below.
\begin{equation}\label{grd25} 
\begin{diagram}
\node{}
\arrow{e,t}{\alpha_\beta}\node{}\arrow{e,t}{\gamma_\delta}\node{}\\
\node{w}\arrow{n,l}{\alpha_\gamma}\arrow{e,b}{\alpha_\beta}\node{}\arrow{n,r}{\beta_\gamma}\arrow{e,b}{\gamma_\delta}\node{}\arrow{n,r}{\beta_\delta}
\end{diagram}
\end{equation}
\end{lemma}

\begin{proof} We check that we cannot have $\beta'=\gamma'$ if the first jeu de taquin rule applied is one of (J0), (J1), (J2), (J5), or (J5$'$). All cases except (J2) are straightforward, involving the application of a rule (J0). In the case of (J2), we have  $\overline{\Delta}=\alpha_x$. The subcase $x\ne\gamma$ is straightforward (the second rule applied is (J0)), while if $x=\gamma$ we necessarily have the growth diagram (\ref{grd25}), because the second rule applied has to be (J5). Indeed, it cannot be (J5$'$), because $\beta$ and $\gamma$ are in positions 1 through $k$ of $(\alpha,\beta)w$ by (J2), where $w$ is the permutation in the bottom left corner of the growth diagram.

This leaves us with the cases when the first jeu de taquin rule applied is (J7$'$) or (J8). In these cases we have $\overline{\Delta}=\beta_x$. In the first case we have $\overline{\Gamma}'=\beta_x\gamma_\delta$, so we obtain the first growth diagram in (\ref{2grd}) if $x=\gamma$. In the second case, we clearly cannot have $\beta'=\gamma'$ if $x\ne \gamma$ (the second jeu de taquin rule applied is (J0)); but if $x=\gamma$ we obtain the second growth diagram, because we cannot use (J5$'$) as the second rule (it leads to $\beta'=\delta'$, which is forbidden), so we have to use (J5).
\end{proof}
 
We now consider the case when $\overline{\Gamma}$ is the left-hand side of a relation (KB3).

\begin{lemma}\label{kb30} Assume that $\overline{\Gamma}$ is the left-hand side of a relation {\rm (KB3)}, that is, $\overline{\Gamma}=\alpha_\beta\varepsilon_\varphi\gamma_\delta$, where $\alpha<\beta\le\gamma<\delta\le\varepsilon<\varphi$. Let $\overline{\Gamma}'=\alpha'_{\beta'}\varepsilon'_{\varphi'}\gamma'_{\delta'}$. Then the only cases when $\beta'=\varepsilon'$ are the ones shown below, in which the bottom rows are the left-hand sides of relations {\rm (KB4)} and {\rm (KB1)}, respectively.
\begin{equation}\label{2grd3} 
\begin{diagram}
\node{}
\arrow{e,t}{\alpha_\beta}\node{}\arrow{e,t}{\varepsilon_\varphi}\node{}\arrow{e,t}{\beta_\varepsilon}\node{}\\
\node{w}\arrow{n,l}{\beta_\varepsilon}\arrow{e,b}{\beta_\varepsilon}\node{}\arrow{n,r}{\alpha_\beta}\arrow{e,b}{\varepsilon_\varphi}\node{}\arrow{n,r}{\alpha_\beta}\arrow{e,b}{\alpha_\beta}\node{}\arrow{n,r}{\beta_\varepsilon}
\end{diagram}\;\;\;\;\;\;\;\;\;\;\;\;\;\;\;\;
\begin{diagram}
\node{}
\arrow{e,t}{\alpha_\beta}\node{}\arrow{e,t}{\varepsilon_\varphi}\node{}\arrow{e,t}{\beta_\varepsilon}\node{}\\
\node{w}\arrow{n,l}{\beta_\varepsilon}\arrow{e,b}{\alpha_\varepsilon}\node{}\arrow{n,r}{\beta_\varepsilon}\arrow{e,b}{\varepsilon_\varphi}\node{}\arrow{n,r}{\beta_\varphi}\arrow{e,b}{\beta_\varepsilon}\node{}\arrow{n,r}{\varepsilon_\varphi}
\end{diagram}
\end{equation}
Moreover, if $\overline{\Gamma}\sim\overline{\Gamma}_r$ and $\overline{\Gamma}'\sim\overline{\Gamma}_r'$, then $(\overline{\Delta},\overline{\Gamma}_r)\mapsto(\overline{\Gamma}_r',\overline{\Delta}')$.
\end{lemma}

\begin{proof}
We will show that $\beta'=\varepsilon'$ implies $\gamma=\beta$ and $\delta=\varepsilon$. Then Lemma \ref{lesstoeq} immediately leads to the two growth diagrams in (\ref{2grd3}). Moreover, replacing $\overline{\Gamma}$ by $\overline{\Gamma}_r$ in the mentioned diagrams, we obtain the corresponding diagrams below; this proves the second part of the lemma.
\begin{equation*}
\begin{diagram}
\node{}
\arrow{e,t}{\varepsilon_\varphi}\node{}\arrow{e,t}{\alpha_\beta}\node{}\arrow{e,t}{\beta_\varepsilon}\node{}\\
\node{w}\arrow{n,l}{\beta_\varepsilon}\arrow{e,b}{\beta_\varepsilon}\node{}\arrow{n,r}{\varepsilon_\varphi}\arrow{e,b}{\alpha_\beta}\node{}\arrow{n,r}{\varepsilon_\varphi}\arrow{e,b}{\varepsilon_\varphi}\node{}\arrow{n,r}{\beta_\varepsilon}
\end{diagram}\;\;\;\;\;\;\;\;\;\;\;\;\;\;\;\;
\begin{diagram}
\node{}
\arrow{e,t}{\varepsilon_\varphi}\node{}\arrow{e,t}{\alpha_\beta}\node{}\arrow{e,t}{\beta_\varepsilon}\node{}\\
\node{w}\arrow{n,l}{\beta_\varepsilon}\arrow{e,b}{\beta_\varepsilon}\node{}\arrow{n,r}{\varepsilon_\varphi}\arrow{e,b}{\alpha_\beta}\node{}\arrow{n,r}{\varepsilon_\varphi}\arrow{e,b}{\beta_\varphi}\node{}\arrow{n,r}{\varepsilon_\varphi}
\end{diagram}
\end{equation*}

We now show that, under the assumption $\beta'=\varepsilon'$, we can rule out the cases when at least one weak inequality in $\alpha<\beta\le\gamma<\delta\le\varepsilon<\varphi$ is a strict inequality. Let $w$ be, as above, the permutation in the bottom left corner of the growth diagram. As above, by Lemma \ref{lesstoeq}, we know that $\overline{\Delta}=\beta_\varepsilon$. Assume first that $\delta<\varepsilon$. By the $k$-semi-shuffle condition, $(\beta,\varepsilon)w$ has the form $\ldots|\ldots\beta\ldots\delta\ldots\varphi\ldots$. But then $w$ is not a $k$-semi-shuffle. So we are left with the case $\beta<\gamma$ and $\delta=\varepsilon$. Now $\beta$ and $\gamma$ are in positions $1$ through $k$ of $w$, while $\varepsilon$ is in a position greater than $k$. The entry $\beta$ cannot precede $\gamma$ in $w$ because the transposition $(\beta,\varepsilon)$ would violate the cover condition. It means that $\gamma$ precedes $\beta$, which leads to a contradiction by Lemma \ref{movevert}.
\end{proof}

We can now prove  Lemma \ref{presplac} when $\overline{\Gamma}$ is the left-hand side of a relation (KB3).

\begin{lemma}\label{kb3}
Let $(\overline{\Delta},\overline{\Gamma})\mapsto (\overline{\Gamma}',\overline{\Delta}')$ be a height $1$ fragment of a growth diagram for which $\overline{\Delta}$  has property {\rm (PL)} and $\overline{\Gamma}$ is the left-hand side of a relation {\rm (KB3)}. Then $\overline{\Gamma}'$ is the left-hand side of a relation {\rm (KB3)}, {\rm (KB4)}, or {\rm (KB1)}. Moreover, if $\overline{\Gamma}\sim\overline{\Gamma}_r$ and $\overline{\Gamma}'\sim\overline{\Gamma}_r'$, then $(\overline{\Delta},\overline{\Gamma}_r)\mapsto(\overline{\Gamma}_r',\overline{\Delta}')$.
\end{lemma}

\begin{proof}
As usual, we let $\overline{\Gamma}=\alpha_\beta\varepsilon_\varphi\gamma_\delta$, where $\alpha<\beta\le\gamma<\delta\le\varepsilon<\varphi$, and $\overline{\Gamma}'=\alpha'_{\beta'}\varepsilon'_{\varphi'}\gamma'_{\delta'}$. By Lemma \ref{kb30}, it suffices to consider the case when $\beta'<\varepsilon'$. By Lemma \ref{gtr}, we have $\gamma'<\delta'\le\varepsilon'<\varphi'$. Let $\overline{\Gamma}_r':=\varepsilon'_{\varphi'}\alpha'_{\beta'}\gamma'_{\delta'}$. Consider first the case when  at least one of the first two jeu de taquin rules applied in the diagram $(\overline{\Delta},\overline{\Gamma})\mapsto (\overline{\Gamma}',\overline{\Delta}')$ is (J0). Then, it is easy to see that $(\overline{\Delta},\overline{\Gamma}_r)\mapsto (\overline{\Gamma}_r',\overline{\Delta}')$. Indeed, if $\overline{\Delta}=x_y$ and the first rule used is (J0), then $\{\alpha,\beta\}\cap\{x,y,\varepsilon,\varphi\}=\emptyset$, and we have a similar property if the second rule is (J0). By Lemma \ref{less}, we now have $\alpha'<\beta'\le\gamma'<\delta'$, so $\overline{\Gamma}'$ is the left-hand side of a relation (KB3) and $\overline{\Gamma}'\sim\overline{\Gamma}_r'$, as sought. 

It remains to investigate the case when none of the first two jeu de taquin rules applied in the diagram $(\overline{\Delta},\overline{\Gamma})\mapsto (\overline{\Gamma}',\overline{\Delta}')$ is (J0). By the second part of Lemma \ref{lesstoeq}, the corresponding growth diagram must be the one below, where $w$ has the form $\ldots\alpha\ldots\varepsilon\ldots|\ldots\beta\ldots\varphi\ldots$.
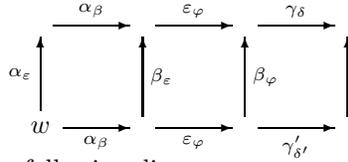
\begin{equation*} 
\begin{diagram}
\node{}
\arrow{e,t}{\alpha_\beta}\node{}\arrow{e,t}{\varepsilon_\varphi}\node{}\arrow{e,t}{\gamma_\delta}\node{}\\
\node{w}\arrow{n,l}{\alpha_\varepsilon}\arrow{e,b}{\alpha_\beta}\node{}\arrow{n,r}{\beta_\varepsilon}\arrow{e,b}{\varepsilon_\varphi}\node{}\arrow{n,r}{\beta_\varphi}\arrow{e,b}{\gamma'_{\delta'}}\node{}\arrow{n,r}{}
\end{diagram}
\end{equation*}
But then we also have the following diagram.
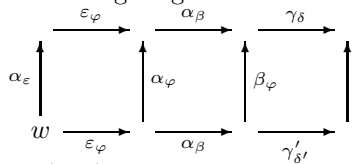
\begin{equation*} 
\begin{diagram}
\node{}
\arrow{e,t}{\varepsilon_\varphi}\node{}\arrow{e,t}{\alpha_\beta}\node{}\arrow{e,t}{\gamma_\delta}\node{}\\
\node{w}\arrow{n,l}{\alpha_\varepsilon}\arrow{e,b}{\varepsilon_\varphi}\node{}\arrow{n,r}{\alpha_\varphi}\arrow{e,b}{\alpha_\beta}\node{}\arrow{n,r}{\beta_\varphi}\arrow{e,b}{\gamma'_{\delta'}}\node{}\arrow{n,r}{}
\end{diagram}
\end{equation*}
Hence, we have $(\overline{\Delta},\overline{\Gamma}_r)\mapsto (\overline{\Gamma}_r',\overline{\Delta}')$ as in the above case, and the same reasoning as before concludes the proof.
\end{proof}

\subsection{The relation (KB4)}

The proof of Lemma \ref{presplac} when $\overline{\Gamma}$ is the left-hand side of a relation (KB4) is completely similar to the proof for the relation (KB3), and therefore is omitted. Nevertheless, let us mention that the new proof is based on the two lemmas below, which are the analogs of Lemmas \ref{lesstoeq} and \ref{kb3}, respectively.

\begin{lemma}
Let $\overline{\Gamma}=\gamma_\delta\alpha_\beta$ and $\overline{\Gamma}'=\gamma'_{\delta'}\alpha'_{\beta'}$, where $\beta<\gamma$. The only case when $\beta'=\gamma'$ is the one shown in the diagram below.
\begin{equation}
\begin{diagram}
\node{}
\arrow{e,t}{\gamma_\delta}\node{}\arrow{e,t}{\alpha_\beta}\node{}\\
\node{}\arrow{n,l}{\beta_\gamma}\arrow{e,b}{\beta_\gamma}\node{}\arrow{n,r}{\gamma_\delta}\arrow{e,b}{\alpha_\beta}\node{}\arrow{n,r}{\gamma_\delta}
\end{diagram}
\end{equation}
If $\beta'<\gamma'$, then there is only one case when none of the two jeu de taquin rules is {\rm (J0)}, namely the growth diagram below, which is the pair of {\rm (\ref{grd25})}.
\begin{equation} 
\begin{diagram}
\node{}
\arrow{e,t}{\gamma_\delta}\node{}\arrow{e,t}{\alpha_\beta}\node{}\\
\node{}\arrow{n,l}{\alpha_\gamma}\arrow{e,b}{\gamma_\delta}\node{}\arrow{n,r}{\alpha_\delta}\arrow{e,b}{\alpha_\beta}\node{}\arrow{n,r}{\beta_\delta}
\end{diagram}
\end{equation}
\end{lemma}

\begin{lemma}\label{kb4}
Let $(\overline{\Delta},\overline{\Gamma})\mapsto (\overline{\Gamma}',\overline{\Delta}')$ be a height $1$ fragment of a growth diagram for which $\overline{\Delta}$  has property {\rm (PL)} and $\overline{\Gamma}$ is the left-hand side of a relation {\rm (KB4)}. Then $\overline{\Gamma}'$ is always the left-hand side of a relation {\rm (KB4)}. Moreover, if $\overline{\Gamma}\sim\overline{\Gamma}_r$ and $\overline{\Gamma}'\sim\overline{\Gamma}_r'$, then $(\overline{\Delta},\overline{\Gamma}_r)\mapsto(\overline{\Gamma}_r',\overline{\Delta}')$.
\end{lemma}

\subsection{The relation (KB1)}

We now prove Lemma \ref{presplac} when $\overline{\Gamma}$ is the left-hand side of a relation (KB1), namely $\overline{\Gamma}=\alpha_\gamma\gamma_\delta\beta_\gamma$ with $\alpha<\beta<\gamma<\delta$. Let $\overline{\Delta}=x_y$, and let $w$ be the permutation in the top left corner of the diagram. By Lemma \ref{ruleoutkb}, $w$ has the form $\ldots\beta\ldots\alpha\ldots|\ldots\gamma\ldots\delta\ldots$. 

We have a nontrivial case only when $\{x,y\}\cap\{\alpha,\beta,\gamma,\delta\}$ is nonempty. We will consider separately the cases when this intersection contains $\alpha$, $\beta$, $\gamma$, and $\delta$. An important observation is that, if $x\in\{\alpha,\beta,\gamma,\delta\}$ and $y\not\in\{\alpha,\beta,\gamma,\delta\}$, then $y$ has to precede $\beta$ in $w$, due to Lemma \ref{movevert}; so $w$ has the form $\ldots y\ldots\beta\ldots\alpha\ldots|\ldots\gamma\ldots\delta\ldots$. Below we represent growth diagrams in pairs: $(\overline{\Delta},\overline{\Gamma})\mapsto (\overline{\Gamma}',\overline{\Delta}')$ and $(\overline{\Delta},\overline{\Gamma}_r)\mapsto(\overline{\Gamma}_r',\overline{\Delta}')$.

{\em Case} 1: $x$ or $y$ is $\alpha$. Consider first the case $x=\alpha$, with $y\not\in\{\alpha,\beta,\gamma,\delta\}$. We must have $\alpha<y<\beta$, because $y>\beta$ would violate the cover condition. This leads to the following growth diagrams.
\begin{equation*}
\begin{diagram}
\node{w}
\arrow{e,t}{\alpha_\gamma}\node{}\arrow{e,t}{\gamma_\delta}\node{}\arrow{e,t}{\beta_\gamma}\node{}\\
\node{}\arrow{n,l}{\alpha_y}\arrow{e,b}{y_\gamma}\node{}\arrow{n,r}{\alpha_y}\arrow{e,b}{\gamma_\delta}\node{}\arrow{n,r}{\alpha_y}\arrow{e,b}{\beta_\gamma}\node{}\arrow{n,r}{\alpha_y}
\end{diagram}\;\;\;\;\;\;\;\;\;\;\;\;\;\;\;\;
\begin{diagram}
\node{}
\arrow{e,t}{\beta_\gamma}\node{}\arrow{e,t}{\alpha_\beta}\node{}\arrow{e,t}{\beta_\delta}\node{}\\
\node{}\arrow{n,l}{\alpha_y}\arrow{e,b}{\beta_\gamma}\node{}\arrow{n,r}{\alpha_y}\arrow{e,b}{y_\beta}\node{}\arrow{n,r}{\alpha_y}\arrow{e,b}{\beta_\delta}\node{}\arrow{n,r}{\alpha_y}
\end{diagram}
\end{equation*}

We cannot have $y=\alpha$, because it contradicts Lemma \ref{movevert}. So the only other possibility is $\overline{\Delta}=\alpha_\beta$, which leads to the following growth diagrams; the latter case is the only one when $\overline{\Gamma}'$ is not the left-hand side of a relation (KB1), being instead the left-hand side of a relation (KB4).

\begin{equation*}
\begin{diagram}
\node{w}\arrow{e,t}{\alpha_\gamma}\node{}
\arrow{e,t}{{\gamma}_{\delta}}\node{}\arrow{e,t}{{\beta}_{\gamma}}\node{}\\
\node{}\arrow{n,l}{\alpha_\beta}\arrow{e,b}{\beta_\gamma}\node{}\arrow{n,r}{\alpha_{\beta}}\arrow{e,b}{{\gamma}_{\delta}}\node{}\arrow{n,r}{\alpha_{\beta}}\arrow{e,b}{\alpha_{\beta}}\node{}\arrow{n,r}{{\beta}_{\gamma}}
\end{diagram}\;\;\;\;\;\;\;\;\;\;\;\;\;\;\;\;
\begin{diagram}
\node{}\arrow{e,t}{\beta_\gamma}\node{}
\arrow{e,t}{\alpha_\beta}\node{}\arrow{e,t}{{\beta}_{\delta}}\node{}\\
\node{}\arrow{n,l}{\alpha_\beta}\arrow{e,b}{\beta_\gamma}\node{}\arrow{n,r}{\alpha_{\gamma}}\arrow{e,b}{\alpha_{\beta}}\node{}\arrow{n,r}{{\beta}_\gamma}\arrow{e,b}{{\gamma}_{\delta}}\node{}\arrow{n,r}{{\beta}_{\gamma}}
\end{diagram}
\end{equation*}

{\em Case} 2: $x$ or $y$ is $\beta$. One possibility is that $y=\beta$, so $x<\beta$. Assuming $x\ne\alpha$, $w$ has the form $\ldots\beta\ldots x\ldots\alpha\ldots|\ldots\gamma\ldots\delta\ldots$ by the cover condition, so we obtain the following growth diagrams. 
\begin{equation*}
\begin{diagram}
\node{w}
\arrow{e,t}{\alpha_\gamma}\node{}\arrow{e,t}{\gamma_\delta}\node{}\arrow{e,t}{\beta_\gamma}\node{}\\
\node{}\arrow{n,l}{x_\beta}\arrow{e,b}{\alpha_\gamma}\node{}\arrow{n,r}{x_\beta}\arrow{e,b}{\gamma_\delta}\node{}\arrow{n,r}{x_\beta}\arrow{e,b}{\beta_\gamma}\node{}\arrow{n,r}{x_\gamma}
\end{diagram}\;\;\;\;\;\;\;\;\;\;\;\;\;\;\;\;
\begin{diagram}
\node{}
\arrow{e,t}{\beta_\gamma}\node{}\arrow{e,t}{\alpha_\beta}\node{}\arrow{e,t}{\beta_\delta}\node{}\\
\node{}\arrow{n,l}{x_\beta}\arrow{e,b}{\beta_\gamma}\node{}\arrow{n,r}{x_\gamma}\arrow{e,b}{\alpha_\beta}\node{}\arrow{n,r}{x_\gamma}\arrow{e,b}{\beta_\delta}\node{}\arrow{n,r}{x_\gamma}
\end{diagram}
\end{equation*}

The other possibility is $x=\beta$. Below we consider separately the cases $\beta<y<\gamma$ and $y>\gamma$ (in the latter case, $y\ne\delta$). 
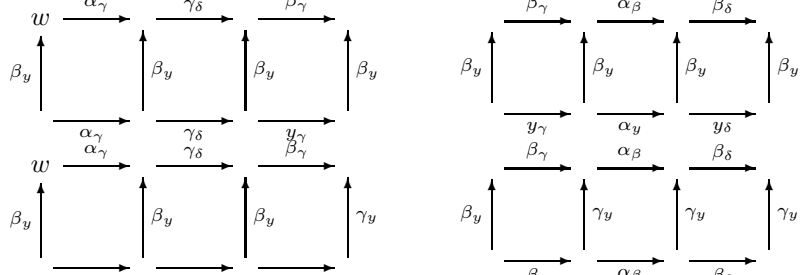
\begin{equation*}
\begin{diagram}
\node{w}
\arrow{e,t}{\alpha_\gamma}\node{}\arrow{e,t}{\gamma_\delta}\node{}\arrow{e,t}{\beta_\gamma}\node{}\\
\node{}\arrow{n,l}{\beta_y}\arrow{e,b}{\alpha_\gamma}\node{}\arrow{n,r}{\beta_y}\arrow{e,b}{\gamma_\delta}\node{}\arrow{n,r}{\beta_y}\arrow{e,b}{y_\gamma}\node{}\arrow{n,r}{\beta_y}
\end{diagram}\;\;\;\;\;\;\;\;\;\;\;\;\;\;\;\;
\begin{diagram}
\node{}
\arrow{e,t}{\beta_\gamma}\node{}\arrow{e,t}{\alpha_\beta}\node{}\arrow{e,t}{\beta_\delta}\node{}\\
\node{}\arrow{n,l}{\beta_y}\arrow{e,b}{y_\gamma}\node{}\arrow{n,r}{\beta_y}\arrow{e,b}{\alpha_y}\node{}\arrow{n,r}{\beta_y}\arrow{e,b}{y_\delta}\node{}\arrow{n,r}{\beta_y}
\end{diagram}
\end{equation*}

\begin{equation*}
\begin{diagram}
\node{w}
\arrow{e,t}{\alpha_\gamma}\node{}\arrow{e,t}{\gamma_\delta}\node{}\arrow{e,t}{\beta_\gamma}\node{}\\
\node{}\arrow{n,l}{\beta_y}\arrow{e,b}{\alpha_\gamma}\node{}\arrow{n,r}{\beta_y}\arrow{e,b}{\gamma_\delta}\node{}\arrow{n,r}{\beta_y}\arrow{e,b}{\beta_\gamma}\node{}\arrow{n,r}{\gamma_y}
\end{diagram}\;\;\;\;\;\;\;\;\;\;\;\;\;\;\;\;
\begin{diagram}
\node{}
\arrow{e,t}{\beta_\gamma}\node{}\arrow{e,t}{\alpha_\beta}\node{}\arrow{e,t}{\beta_\delta}\node{}\\
\node{}\arrow{n,l}{\beta_y}\arrow{e,b}{\beta_\gamma}\node{}\arrow{n,r}{\gamma_y}\arrow{e,b}{\alpha_\beta}\node{}\arrow{n,r}{\gamma_y}\arrow{e,b}{\beta_\delta}\node{}\arrow{n,r}{\gamma_y}
\end{diagram}
\end{equation*}

{\em Case} 3: $x$ or $y$ is $\gamma$. Here we are forced to have $x=\gamma$ and $\gamma<y<\delta$. This leads to the following growth diagrams.
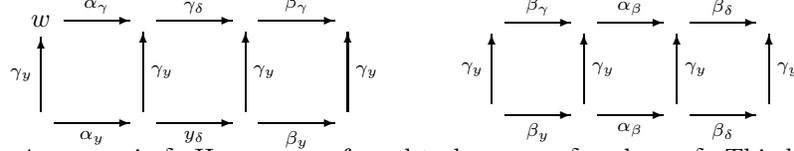
\begin{equation*}
\begin{diagram}
\node{w}
\arrow{e,t}{\alpha_\gamma}\node{}\arrow{e,t}{\gamma_\delta}\node{}\arrow{e,t}{\beta_\gamma}\node{}\\
\node{}\arrow{n,l}{\gamma_y}\arrow{e,b}{\alpha_y}\node{}\arrow{n,r}{\gamma_y}\arrow{e,b}{y_\delta}\node{}\arrow{n,r}{\gamma_y}\arrow{e,b}{\beta_y}\node{}\arrow{n,r}{\gamma_y}
\end{diagram}\;\;\;\;\;\;\;\;\;\;\;\;\;\;\;\;
\begin{diagram}
\node{}
\arrow{e,t}{\beta_\gamma}\node{}\arrow{e,t}{\alpha_\beta}\node{}\arrow{e,t}{\beta_\delta}\node{}\\
\node{}\arrow{n,l}{\gamma_y}\arrow{e,b}{\beta_y}\node{}\arrow{n,r}{\gamma_y}\arrow{e,b}{\alpha_\beta}\node{}\arrow{n,r}{\gamma_y}\arrow{e,b}{\beta_\delta}\node{}\arrow{n,r}{\gamma_y}
\end{diagram}
\end{equation*}

{\em Case} 4: $x$ or $y$ is $\delta$. Here we are forced to have $x=\delta$ and $y>\delta$. This leads to the following growth diagrams.
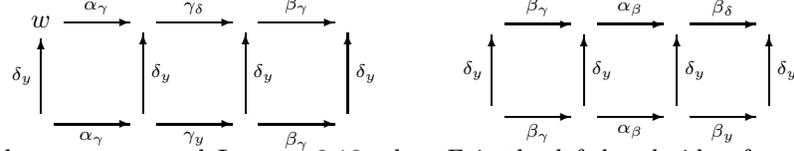
\begin{equation*}
\begin{diagram}
\node{w}
\arrow{e,t}{\alpha_\gamma}\node{}\arrow{e,t}{\gamma_\delta}\node{}\arrow{e,t}{\beta_\gamma}\node{}\\
\node{}\arrow{n,l}{\delta_y}\arrow{e,b}{\alpha_\gamma}\node{}\arrow{n,r}{\delta_y}\arrow{e,b}{\gamma_y}\node{}\arrow{n,r}{\delta_y}\arrow{e,b}{\beta_\gamma}\node{}\arrow{n,r}{\delta_y}
\end{diagram}\;\;\;\;\;\;\;\;\;\;\;\;\;\;\;\;
\begin{diagram}
\node{}
\arrow{e,t}{\beta_\gamma}\node{}\arrow{e,t}{\alpha_\beta}\node{}\arrow{e,t}{\beta_\delta}\node{}\\
\node{}\arrow{n,l}{\delta_y}\arrow{e,b}{\beta_\gamma}\node{}\arrow{n,r}{\delta_y}\arrow{e,b}{\alpha_\beta}\node{}\arrow{n,r}{\delta_y}\arrow{e,b}{\beta_y}\node{}\arrow{n,r}{\delta_y}
\end{diagram}
\end{equation*}

We have now proved Lemma \ref{presplac} when $\Gamma$ is the left-hand side of a relation (KB1). We make a more precise statement below. 

\begin{lemma}\label{kb1}
Let $(\overline{\Delta},\overline{\Gamma})\mapsto (\overline{\Gamma}',\overline{\Delta}')$ be a height $1$ fragment of a growth diagram for which $\overline{\Delta}$  has property {\rm (PL)} and $\overline{\Gamma}$ is the left-hand side of a relation {\rm (KB1)}. Then $\overline{\Gamma}'$ is the left-hand side of a relation {\rm (KB1)} or {\rm (KB4)}. Moreover, if $\overline{\Gamma}\sim\overline{\Gamma}_r$ and $\overline{\Gamma}'\sim\overline{\Gamma}_r'$, then $(\overline{\Delta},\overline{\Gamma}_r)\mapsto(\overline{\Gamma}_r',\overline{\Delta}')$.
\end{lemma}

\subsection{The property (PR)} It is possible to reduce the part of Lemma \ref{presplac} related to property (PR) to the one related to (PL). The idea is to use the automorphism  $w\mapsto w_\circ ww_\circ$ of the Bruhat order on $S_n$, where $w_\circ$ is the longest permutation in $S_n$. Notice that this automorphism interchanges
\begin{itemize}
\item chains in $k$-Bruhat order with chains in $(n-k)$-Bruhat order,
\item $k$-semi-shuffles with permutations having no descents before $n-k$,
\item chains having properties (PL) and (PR),
\item one side of (KB1) with the opposite side of (KB2), as well as the two sides of (KB3) and (KB4) among themselves,
\item the jeu de taquin relations as follows: J0 with itself, J1 with J3, J2 with J4, J5 with J7, J5$'$ with J7$'$, J6 with J8, and J6$'$ with J8$'$. 
\end{itemize}

\bibliography{copies1}

\end{document}